\documentclass[a4paper,12pt]{lms}
\usepackage{amsfonts}
\usepackage{amssymb}
\usepackage{amsmath}
\usepackage{bbm}
\usepackage{epsfig}

\newtheorem{theorem}{Theorem}
\newtheorem{proposition}{Proposition}

\newtheorem{remark}{Remark}
\newtheorem{corollary}{Corollary}
\newtheorem{lemma}{Lemma}

\newcommand{\R}{\Bbb R}

\newcommand{\1}{\mathbbm 1}

\author{Artur Nicolau, Daniel Seco}
\title{Smoothness of sets in Euclidean spaces}
\date{}

\begin{document}
\maketitle


\footnotetext{Both authors are supported partially by the grants
MTM2008-00145 and 2009SGR420.}

\begin{abstract}
We study some properties of smooth sets in the sense defined by
Hungerford. We prove a sharp form of Hungerford's Theorem on the
Hausdorff dimension of their boundaries on Euclidean spaces and show
the invariance of the definition under a class of automorphisms of
the ambient space.
\end{abstract}

\section{Introduction}

$\quad$ The Lebesgue Density Theorem tells us that the density of a
measurable set approximates the characteristic function of the set
at almost every point. We are going to study sets whose densities at
small scales vary uniformly.

In this paper, a cube will mean a cube in the Euclidean space
$\R^{n}$ with sides parallel to the axis. Two cubes $Q, Q' \subset
\R^{n}$ with the same sidelength $l(Q)=l(Q')$ are called consecutive
if the intersection of their closures is one of their faces. Given a
measurable set $A \subset \R^{n}$, let $|A|$ denote its Lebesgue
measure and $D(Q)$ its density in a cube $Q \subset \R^{n}$, that
is, $D(Q)= |A \cap Q| / |Q|$. A measurable set $A \subset \R^{n}$ is
called \emph{smooth} (in $\R^{n}$) if \[\lim_{\delta \rightarrow 0}
\sup |D(Q) - D(Q')| = 0\] where the supremum is taken over all pairs
of consecutive cubes $Q, Q'$ with $l(Q)=l(Q') \leq \delta$. In
dimension $n=1$, this notion was introduced by Hungerford (\cite{2})
in relation to the small Zygmund class. Actually, a set $A \subset
\R$ is smooth if and only if its distribution function is in the
small Zygmund class, or equivalently, the restriction of the
Lebesgue measure to the set $A$ is a smooth measure in the sense of
Kahane (\cite{4}).

Sets $A \subset \R^{n}$ with $|A|=0$ or $|\R^{n}\backslash A|=0$ are
trivially smooth but Hungerford provided non trivial examples, using
a nice previous recursive construction by Kahane. See \cite{2}.
Other sharper examples are given in \cite{1}.

In dimension $n=1$, Hungerford proved that the boundary of a
nontrivial smooth set has full Hausdorff dimension (\cite{2}, see
also \cite{6}). His argument shows that if $A$ is a smooth set in
$\R$ with $|A|> 0$ and $|\R \backslash A| > 0$, then the set of
points $x \in \R$ for which there exists a sequence of intervals
$\{I_{j}\}$ containing $x$ such that \[\lim_{j \rightarrow \infty}
D(I_{j}) = 1/2\] still has Hausdorff dimension $1$. The main goal of
this paper is to sharpen this result and to extend it to Euclidean
spaces. It is worth mentioning that Hungerford arguments cannot be
extended to several dimensions since it is used that the image under
a nontrivial linear mapping of an interval is still an interval, or
more generally, that an open connected set is an interval and this
obviously does not hold for cubes in $\R^{n}$, for $n>1$. Given a
point $x \in \R^{n}$ and $h>0$, let $Q(x, h)$ denote the cube
centered at $x$ of sidelength $h$. With this notation, our main
result is the following:

\begin{theorem}\label{th2}
Let $A$ be a smooth set in $\R^{n}$ with $|A|>0$ and
$|\R^{n}\backslash A| >0$. Fix $0<\alpha<1$. Then the set \[E(A,
\alpha)= \left\{x \in \R^{n}: \lim_{h \rightarrow 0}
D(Q(x,h))=\alpha\right\}\] has Hausdorff dimension $n$.
\end{theorem}

Our result is local, meaning, given a cube $Q \subset \R^{n}$ with
$0 < |A \cap Q| < |Q|$, then $E(A, \alpha)\cap Q$ has full Hausdorff
dimension. As a consequence the Hausdorff dimension of $\partial A
\cap Q$ is $n$.

Section 2 contains a proof of Theorem $\ref{th2}$. A Cantor type
subset of $E(A, \alpha)$ will be constructed and its dimension will
be computed using a standard result. The generations of the Cantor
set will be defined recursively by means of a stopping time
argument. The good averaging properties of the density are used to
estimate the dimension of the Cantor set.

The definition of smooth set concerns the behavior of the density of
the set on the grid of cubes in $\R^{n}$ with sides parallel to the
axis. We consider two natural questions arising from this
definition. First, we study how the definition depends on the grid
of cubes, that is, if other natural grids, such as dyadic cubes or
general parallelepipeds would lead to the same notion. Second, we
consider whether the class of smooth sets is preserved by regular
mappings. It turns out that these questions are related and in
Section 3 we provide a positive answer to both of them. A mapping
$\phi: \R^{n} \rightarrow \R^{n}$ is bilipschitz if there exists a
constant $C \geq 1$ such that $C^{-1}||x-y|| \leq ||\phi(x) -
\phi(y)|| \leq C||x-y||$ for any $x, y \in \R^{n}$.

\begin{theorem}\label{th3}
Let $\phi : \R^{n} \rightarrow \R^{n}$ be a bilipschitz
$\mathcal{C}^{1}$ mapping with uniformly continuous Jacobian. Let $A
\subset \R^{n}$ be a measurable set. Then the following are
equivalent:
\begin{description}
\item [(a)] $A$ is a smooth set
\item [(b)] $\phi^{-1}(A)$ is a smooth set
\item [(c)] $A$ verifies the smoothness condition taking, instead of
the grid of cubes, their images through $\phi$, that is, \[\lim_{|Q|
\rightarrow 0} \frac{|A \cap \phi(Q)|}{|\phi(Q)|} - \frac{|A \cap
\phi (Q')|}{|\phi(Q')|} = 0\]
\end{description}
\end{theorem}

As part (c) states, one could replace in the definition of smooth
set, the grid of cubes by other grids such as the grid of dyadic
cubes or the grid of general parallelepipeds with bounded
eccentricity whose sides are not necessarily parallel to the axis or
even the pullback by $\phi$ of the grid of cubes. One can combine
Theorems $\ref{th2}$ and $\ref{th3}$ to conclude the following:

\begin{corollary}
Let $A$ be a smooth set in $\R^{n}$ with $|A| > 0$ and $|\R^{n}
\backslash A| >0$. Let $\phi: \R^{n} \rightarrow \R^{n}$ be a
bilipschitz $\mathcal{C}^{1}$ mapping with uniformly continuous
jacobian. Fix $0< \alpha < 1$. Then the set \[\left\{x\in \R^{n} :
\lim_{h \rightarrow 0} \frac{|A \cap \phi(Q(x,h))|}{|\phi(Q(x,h))|}=
\alpha \right\}\] has Hausdorff dimension $n$.
\end{corollary}

\textbf{Acknowledgements.} We would like to thank Professors F.
Ced{\'o} and W. Dicks for their help with linear algebra techniques.

\section{Proof of Theorem $\ref{th2}$}

\subsection{Preliminary results}

We begin with a preliminary result on the Hausdorff dimension of
certain Cantor type sets which will be used in the proof of Theorem
$\ref{th2}$. In dimension $n=1$, the result was given by Hungerford
in \cite{2}. See also Theorem 10.5 in \cite{7}. The proof in the
higher dimensional case only requires minor adjustments and it will
be omitted.

\begin{lemma}\label{lem1} For $s=0,1,2,...$ let $G(s)$ be a
collection of closed dyadic cubes in $\R^{n}$ with pairwise disjoint
interiors. Assume that the families are nested, that is \[\bigcup_{Q
\in G(s+1)} Q \subseteq \bigcup_{Q \in G(s)} Q\] Suppose that there
exist two positive constants $0 < P < C < 1$ such that the following
two conditions hold:
\begin{description}
\item [(a)] For any cube $Q \in G(s+1)$ with $Q \subset
\widetilde{Q} \in G(s)$ one has $|Q| \leq P|\widetilde{Q}|$.
\item [(b)] For any $\widetilde{Q} \in G(s)$ one has \[\sum |Q| \geq
C|\widetilde{Q}|\] where the sum is taken over all cubes $Q \in
G(s+1)$ contained in $\widetilde{Q}$.
\end{description}
Let $E(s) = \bigcup Q$, where the union is taken over all cubes in
$G(s)$ and $E \equiv \bigcap_{s=0}^{\infty} E(s)$. Then $\dim E \geq
n (1-\log_{P}C)$.
\end{lemma}

The next auxiliary result is the building block of the Cantor set on
which the set has a fixed density. Recall that given a measurable
set $A \subset \R^{n}$, its density on a cube $Q$ is $D(Q)=|A \cap
Q| /|Q|$. Given a continuous increasing function $\omega: [0,1]
\rightarrow [0, \infty)$ with $\omega(0)=0$, a set $A \subset
\R^{n}$ is called \emph{$\omega$-smooth} (in $\R^{n}$) if \[|D(Q) -
D(Q')| \leq \omega(l(Q))\] for any pair of consecutive cubes $Q,
Q'\subset \R^{n}$ of sidelength $l(Q)=l(Q')$.

\begin{lemma}\label{lem2}
Let $A$ be an $\omega$-smooth set of $\R^{n}$ with $0 < |A \cap
[0,1]^{n}| < 1$. Let $Q$ be a dyadic cube. Fix a constant
$\varepsilon > 0$ such that $n\omega(l(Q)) < \varepsilon < \min
\{D(Q), 1-D(Q)\}$. Let $\mathfrak{A}(Q)$ be the family of maximal
dyadic cubes $Q_{k}$ contained in $Q$ such that
\begin{equation}\label{eq31} |D(Q_{k})- D(Q)| \geq \varepsilon
\end{equation} Then:
\begin{description}
\item [(a)] For any
$Q_{k} \in \mathfrak{A}(Q)$ one has \[|Q_{k}| \leq 2^{- \varepsilon
/ \omega(l(Q))}|Q|\]
\item [(b)] Let $\mathfrak{A}^{+} (Q)$ (respectively
$\mathfrak{A}^{-}(Q)$) be the subfamily of $\mathfrak{A}(Q)$ formed
by those cubes $Q_{k} \in \mathfrak{A}(Q)$ for which $D(Q_{k})-D(Q)
\geq \varepsilon$ (respectively $D(Q)- D(Q_{k}) \geq \varepsilon$).
Then \[\sum |Q_{k}| \geq |Q|/4\] where the sum is taken over all the
cubes $Q_{k} \in \mathfrak{A}^{+}(Q)$ (respectively $Q_{k} \in
\mathfrak{A}^{-}(Q)$).
\end{description}
\end{lemma}

\begin{proof}
If $Q_{1} \subset Q_{2} \subset Q$ are two dyadic cubes with
$l(Q_{1})= l(Q_{2})/2$ then $|D(Q_{1}) - D(Q_{2})| \leq n \omega
(l(Q))$. So if $\eqref{eq31}$ holds, one deduces that
$\log_{2}l(Q_{k})^{-1} \geq \log_{2}l(Q)^{-1} + \varepsilon/n
\omega(l(Q))$. Therefore, (a) is proved.

To prove (b), we observe first that, by Lebesgue Density Theorem,
one has $\displaystyle{\sum_{\mathfrak{A}(Q)} |Q_{k}| = |Q|}$. Also,
\begin{equation}\label{eq32}
\sum_{\mathfrak{A}(Q)} (D(Q_{k}) - D(Q))|Q_{k}| = 0
\end{equation}

We argue by contradiction. Assume that $\displaystyle{\sum_{Q_{k}
\in \mathfrak{A}^{+}(Q)} |Q_{k}|} < |Q|/4$ and hence
$\displaystyle{\sum_{Q_{k} \in \mathfrak{A}^{-}(Q)} |Q_{k}|} \geq
3|Q|/4$, which gives us \[\sum_{Q_{k} \in \mathfrak{A}^{-}(Q)}
(D(Q_{k}) - D(Q))|Q_{k}| \leq -3\varepsilon|Q|/4\] The maximality of
$Q_{k}$ tells us that $|D(Q_{k})-D(Q)| \leq \varepsilon + n
\omega(l(Q)) < 2 \varepsilon$. Therefore \[\sum_{Q_{k} \in
\mathfrak{A}^{+}(Q)} (D(Q_{k})-D(Q)) |Q_{k}| \leq \varepsilon|Q|/2\]
which contradicts $\eqref{eq32}$. The same argument works for
$\mathfrak{A}^{-}(Q)$.
\end{proof}

\subsection{The dyadic case}

Our next goal is to prove a dyadic version of Theorem $\ref{th2}$,
which already contains its core. Let $Q_{k}(x)$ be the dyadic cube
of generation $k$ which contains the point $x \in \R^{n}$.

\begin{proposition}\label{dyadsmooth}
Let $A$ be a smooth set in $\R^{n}$ with $0 < |A \cap [0,1]^{n}| <
1$. For $0 < \alpha < 1$ consider the set $E_{1}(A, \alpha) =
\left\{x \in [0,1]^{n} : \displaystyle{\lim_{k \rightarrow \infty}}
D(Q_{k}(x)) = \alpha \right\}$. Then $\dim E_{1}(A, \alpha) = n$.
\end{proposition}
\begin{proof}
Fix $0 < \alpha < 1$. A Cantor type set contained in $E_{1}(A,
\alpha)$ will be constructed and Lemma \ref{lem1} will be used to
calculate its Hausdorff dimension. The Cantor type set will be
constructed using generations $G(s)$ which will be defined using
Lemma $\ref{lem2}$, yielding the estimates appearing in Lemma
$\ref{lem1}$.

Given the smooth set $A \subset \R^{n}$ consider the function
\[\omega(t) = \sup |D(Q) - D(Q')|, \quad 0 < t \leq 1\] where the
supremum is taken over all pair of consecutive cubes $Q$ and $Q'$ of
the same sidelength $l(Q) = l(Q') \leq t$. Observe that $\lim
\omega(t)=0$ as $t \rightarrow 0$. Pick a positive integer $k_{0}$
such that $\omega(2^{-k_{0}}) < \min \{\alpha, 1 -\alpha\}/20$.
Define an increasing sequence $\{c_{k}\}$ with $c_{k} \rightarrow
\infty$ as $k \rightarrow \infty$ and $c_{k} \geq 2n$ for any $k$,
satisfying $\varepsilon_{k}=c_{k}\omega(2^{-k-k_{0}}) \rightarrow 0$
as $k \rightarrow \infty$. We can also assume $\varepsilon_{k} <
\min \{\alpha, 1-\alpha\}/10$ for any $k=1,2...$ Since $0< |A \cap
[0,1]^{n}|< 1$, there are some small dyadic cubes in $[0,1]^{n}$
with density close to $0$ and others with density close to $1$.
Since $A$ is smooth, we can choose a dyadic cube $Q_{1}$ with
$l(Q_{1}) \leq 2^{-k_{0}-1}$ and $|D(Q_{1}) - \alpha| <
\varepsilon_{1}/2$. Then define the first generation
$G(1)=\{Q_{1}\}$. The next generations are constructed inductively
as follows. Assume that the $k$-th generation $G(k)$ has been
defined so that the following two conditions are satisfied: $l(Q)
\leq 2^{-k-k_{0}}$ and $|D(Q)- \alpha| < \varepsilon_{k}/2$ for any
cube $Q \in G(k)$. The generation $G(k+1)$ is constructed in two
steps. Roughly speaking, we first find cubes whose density is far
away from $\alpha$ and later we find subcubes with density close to
$\alpha$. For $Q \in G(k)$ consider the family $\mathfrak{R}(Q)$ of
maximal dyadic cubes $R \subset Q$ such that $|D(R) - D(Q)| \geq
\varepsilon_{k}$. Observe that, by Lebesgue Density Theorem, $\sum
|R| = |Q|$, where the sum is taken over all cubes $R \in
\mathfrak{R}(Q)$. Fix $R \in \mathfrak{R}(Q)$. Since the set A is
$\omega$-smooth, the difference of densities between two dyadic
cubes $Q_{1} \subset Q_{2} \subset Q$, with $l(Q_{2})=2l(Q_{1})$, is
smaller than $n \omega (l(Q))$. Hence to achieve such a cube $R$ we
need to go through at least $\varepsilon_{k}/ n
\omega(2^{-k-k_{0}})=c_{k}/n$ dyadic steps. Hence
\begin{equation}\label{eq35} |R| \leq 2^{-c_{k}}|Q|
\end{equation} The maximality and the estimate $l(Q) \leq 2^{-k-k_{0}}$ give that $|D(R)
- D(Q)| \leq \varepsilon_{k} + n \omega(2^{-k-k_{0}-1})$. Since
\[|D(R) - \alpha| > \varepsilon_{k}/2 > n \omega(2^{-k-k_{0}}) \geq
n \omega(l(R))\] one can apply Lemma \ref{lem2} with the parameter
$\varepsilon=|D(R)-\alpha|$. In this way, one obtains two families
$\mathfrak{A}^{-}(R)$ and $\mathfrak{A}^{+}(R)$ of dyadic cubes
contained in $R$, according to whether their densities are smaller
or bigger than $D(R)$, but we will only be interested on one of them
which will be called $G_{k+1}(R)$. If $D(R)> \alpha$ we choose
$G_{k+1}(R)=\mathfrak{A}^{-}(R)$. Otherwise, take $G_{k+1}(R)
=\mathfrak{A}^{+}(R)$. Fix, now, $Q^{*} \in G_{k+1}(R)$. The
maximality gives that $|D(Q^{*}) - \alpha| \leq n \omega (l(R))$.
Since $l(R) \leq 2^{-k-k_{0}-1}$ we deduce that $|D(Q^{*})- \alpha|
< \varepsilon_{k+1}/2$. Also, $l(Q^{*}) < l(R)/2 \leq
2^{-k-k_{0}-1}$. Notice that any dyadic cube $\widetilde{Q}$ with
$Q^{*} \subset \widetilde{Q} \subset Q$ satisfies
\begin{equation}\label{eq33} |D(\widetilde{Q}) - \alpha| \leq 6
\varepsilon_{k}
\end{equation} The generation $G(k+1)$ is defined as
\[G(k+1)= \bigcup_{Q \in G(k)} \bigcup_{R \in \mathfrak{R}(Q)}
G_{k+1}(R)\]

Next we will compute the constants appearing in Lemma $\ref{lem1}$.
Let $Q \in G(k)$ and $R \in \mathfrak{R}(Q)$. Part (b) of Lemma
$\ref{lem2}$ says that $\sum|Q_{j}| \geq |R|/4$ where the sum is
taken over all cubes $Q_{j} \in G_{k+1}(R)$. Since
$\displaystyle{\sum_{\mathfrak{R}(Q)}} |R| = |Q|$ one deduces that
\begin{equation}\label{eq34} \sum |Q_{j}| \geq |Q|/4
\end{equation} where the sum is taken over all cubes $Q_{j} \in G(k+1), Q_{j} \subset Q$.
Also, if $Q_{j} \in G(k+1)$ and $Q_{j} \subset Q \in G(k)$, estimate
$\eqref{eq35}$ guarantees that \begin{equation}\label{eq36} |Q_{j}|
\leq 2^{-c_{k}}|Q|
\end{equation} For $k=1,2...$ let $E(k)$ be the union of the cubes of the family
$G(k)$ and let $E = \bigcap E(k)$ be the corresponding Cantor type
set.

Next we show that $E$ is contained in $E_{1}(A, \alpha)$. To do
this, fix $x \in E$ and for any $k=1,2,...$ pick the cube $Q_{k} \in
G(k)$ containing $x$. Let $Q$ be a dyadic cube which contains $x$,
and $k$ the integer for which $Q_{k+1} \subset Q \subset Q_{k}$.
Observe that $k \rightarrow \infty$ as $l(Q) \rightarrow 0$. By
$\eqref{eq33}$ one deduces that $|D(Q) - \alpha| \leq 6
\varepsilon_{k}$ and therefore $x \in E_{1}(A, \alpha)$.

Finally, we apply Lemma $\ref{lem1}$ to show that the dimension of
$E$ is $n$. Actually $\eqref{eq34}$ and $\eqref{eq36}$ give that one
can take $C=1/4$ and $P=2^{-c_{k}}$ in Lemma \ref{lem1}. Hence, the
dimension of $E$ is bigger than $n(1-(2/c_{k}))$. Since $c_{k}
\rightarrow \infty$ as $k \rightarrow \infty$, we deduce that $\dim
E = n$.
\end{proof}

\subsection{Affine control and Proof of Theorem $\ref{th2}$}

We want to study the density of a smooth set in non dyadic cubes.
The proof of Theorem $\ref{th2}$ will be based on Proposition
$\ref{dyadsmooth}$ and on the following auxiliary result on the
behavior of the densities with respect to affine perturbations.

\begin{lemma} \label{lem3}
Let $A$ be a smooth set in $\R^{n}$. Consider the function
\[\omega(t) = \sup |D(Q)-D(Q')|\] where the supremum is taken over
all pairs of consecutive cubes $Q, Q' \subset \R^{n}$, of sidelength
$l(Q) = l(Q') \leq t$.
\begin{description}
\item[(a)] Let $Q, \widetilde{Q}$ be two cubes in $\R^{n}$ with non
empty intersection such that $l(Q)=l(\widetilde{Q})$. Then $|D(Q) -
D(\widetilde{Q})| \leq 3n^{2} \omega (l(Q))$
\item[(b)] Let $Q$ be a cube in $\R^{n}$ and let $tQ$ denote the cube
with the same center as $Q$ and sidelength $tl(Q)$. Then for any $1
\leq t \leq 2$ one has \[|D(Q) - D(tQ)| \leq c(n)\omega(l(Q))\] Here
$c(n)$ is a constant which only depends on the dimension.
\end{description}
\end{lemma}

\begin{proof}
To prove (a), suppose, without loss of generality, that
$Q=[0,1]^{n}$. Let $Q'=[x,1+x] \times [0,1]^{n-1}$ where $-1 < x <
1$. We will show that \begin{equation} \label{eq37} |D(Q) - D(Q')|
\leq 3 n \omega(1)
\end{equation} Since any $Q'$ intersecting $Q$ is of the form $Q'=[x_{1}, 1+x_{1}]
\times ... \times [x_{n}, 1+x_{n}]$ part (a) follows using
$\eqref{eq37}$ $n$ times.

To show $\eqref{eq37}$, decompose $[x,1+x]$ into dyadic intervals,
that is, $[x, 1+x] = \bigcup I_{k}$, where $I_{k}$ is a dyadic
interval of length $2^{-k}$ for $k=1,2,...$ Consider the
parallelepiped $R_{k} = I_{k} \times [0,1]^{n-1}$ and the density
$D(R_{k})$ of the set $A$ on $R_{k}$, meaning, $D(R_{k}) = |R_{k}
\cap A| / |R_{k}|$, $k=1,2...$. The set $R_{k}$ can be split into a
family $\mathfrak{F}_{k}$ of $2^{k (n-1)}$ pairwise disjoint cubes S
of sidelength $2^{-k}$. Since $|D(S) - D(Q)| \leq n \omega (1)
(k+1)$ for any $S \in \mathfrak{F}_{k}$ and $D(R_{k})$ is the mean
of $D(S), S \in \mathfrak{F}_{k}$, we deduce that $|D(R_{k}) - D(Q)|
\leq n \omega(1) (k+1)$. Since \[D(Q')= \sum_{k=1}^{\infty}
2^{-k}D(R_{k})\] we deduce that $|D(Q') - D(Q)| \leq 3n \omega(1)$
which is $\eqref{eq37}$.

We turn now to (b). We can assume $Q=[-1/2,1/2]^{n}$. Consider the
binary decomposition of $t$, that is, $t=
1+\displaystyle{\sum_{k=1}^{\infty}} t_{k}2^{-k}$, with $t_{k} \in
\{0,1\}$. For $m=1,2,...$ let $Q_{m}$ be the cube with the same
center as $Q$ but with sidelength $l(Q_{m})= 1 +
\displaystyle{\sum_{k=1}^{m}} t_{k}2^{-k}$. Then $tQ=
\displaystyle{\bigcup_{m=0}^{\infty}} R_{m}$ where $R_{0} \equiv
Q_{0} \equiv Q$ and $R_{m}=Q_{m} \backslash Q_{m-1}$ for $m \geq 1$.
So for $m \geq 1, R_{m}$ is empty whenever $t_{m}=0$ and otherwise
we will estimate $|D(R_{m}) - D(Q)|$. With this aim assume that
$t_{m}=1$ and split $R_{m}$ into a family $\mathfrak{F}_{m}$ of
pairwise disjoint cubes S of sidelength $2^{-m-1}$. Thus, $|D(S) -
D(Q)| \leq (n(m+1)+1) \omega (l(Q))$ for any $S \in
\mathfrak{F}_{m}$. Since $D(R_{m})$ is the mean of $D(S), S \in
\mathfrak{F}_{m}$, this implies that
\begin{equation}\label{eq38} |D(R_{m}) - D(Q)| \leq (n(m+1)+1)
\omega (l(Q)) \quad m=1,2...
\end{equation} As we also have \[D(tQ) = \sum_{m=0}^{\infty}
\frac{|R_{m}|}{|tQ|}D(R_{m})\] using $\eqref{eq38}$ we obtain
\[|D(tQ) - D(Q)| \leq \omega (l(Q)) \sum_{m=0}^{\infty}
\frac{|R_{m}|}{|tQ|}(n(m+1)+1)\] Since $|R_{m}| \leq C(n) 2^{-m}$
the sum is convergent and the proof is complete.
\end{proof}

We are now ready to prove Theorem $\ref{th2}$.

\begin{proof}[of Theorem $\ref{th2}$]
Applying Proposition $\ref{dyadsmooth}$, we only need to check that
\[\lim_{h \rightarrow 0} D(Q(x,h)) = \alpha\] for any $x \in E_{1}(A,
\alpha)$. Given $h > 0$, let $k$ be the unique integer such that
$2^{-k} \leq h < 2^{-k+1}$. Consider the cube $h2^{k}Q_{k}(x)$ and
apply Lemma $\ref{lem3}$ to deduce that \[\lim_{h \rightarrow 0}
|D(Q(x,h)) - D(Q_{k}(x))| = 0\]
\end{proof}

\section{Equivalent definitions and invariance}

The definition of a smooth set involves the density of the set on
the grid of all cubes in $\R^{n}$ with sides parallel to the axis.
The main purpose of this section is to study the situation for
perturbations of this grid of cubes. The first step consists of
considering linear deformations of the family of cubes, obtaining a
certain grid of parallelepipeds in $\R^{n}$. Afterwards, we will
consider the more general case of the grid arising from a
bilipschitz image of the family of cubes.

\begin{proposition} \label{prop41}
Let $\phi : \R^{n} \rightarrow \R^{n}$ be a linear mapping, and let
$A \subset \R^{n}$ be a smooth set. Then $\phi (A)$ is smooth.
\end{proposition}

\begin{proof}
We can assume that $\phi$ is a linear isomorphism. Given the smooth
set $A$ consider \[\omega(t) = \sup |D(Q)-D(Q')|, \quad 0 < t \leq
1\] where the supremum is taken over all pairs of consecutive cubes
$Q$ and $Q'$ of the same sidelength $l(Q) = l(Q') \leq t$. Since
$|\phi(A) \cap Q| = c(\phi)|A \cap \phi^{-1}(Q)|$, where $c(\phi)$
is a constant which only depends on $\phi$, it is sufficient to show
the smoothness condition, taking, instead of cubes, their preimages
through $\phi$, that is
\begin{equation}\label{eq41} \lim_{|Q| \rightarrow 0} \frac{|A \cap
\phi^{-1}(Q)|- |A \cap \phi^{-1}(Q')|}{|Q|} = 0
\end{equation} Apply the Singular Value Decomposition (see, for instance, Theorem
7.3.5 on \cite{5}) to $\phi$, which allows us to write $\phi= V
\Sigma W$ where $V, W$ are orthogonal mappings and $\Sigma$ is a
diagonal mapping, that is, the matrix of $\Sigma$ is diagonal.
Moreover the elements of the diagonal of $\Sigma$ are the positive
square roots of the eigenvalues of $\phi \phi^{*}$. It will prove
useful later that the elements $\lambda \in \R$ of the diagonal of
$\Sigma$ verify $||\phi^{-1}|| \leq |\lambda| \leq ||\phi||$, where
$||\phi||$ denotes the norm of $\phi$ as a linear mapping from
$\R^{n}$ to $\R^{n}$. So, we proceed to prove the Proposition for
these two cases: (a) $\phi$ is an orthogonal mapping and (b) $\phi$
is a dilation.

We study first case (a). As any orthogonal application is either a
rotation or the composition of a rotation and a reflection by a
subspace parallel to the axis (leaving invariant the grid of cubes),
we can reduce the orthogonal case to rotations. Furthermore, we can
assume that $\phi$ is the identity on a subspace of dimension $n-2$
generated by elements of the canonical basis of $\R^{n}$ and a
rotation of angle $\alpha \in [\pi/6, \pi/4]$ on its orthogonal
complement. Actually, any rotation can be written as the composition
of at most $3n(n-1)$ rotations of this form (see, for instance,
\cite{5}). Let $\widetilde{Q}$ (respectively, $\widetilde{Q'}$) be
the cube centered at the center of $\phi^{-1}(Q)$ (respectively,
$\phi^{-1}(Q')$) of sidelength $l(Q)$ with sides parallel to the
axis. Observe that $\widetilde{Q'}$ is a translation of
$\widetilde{Q}$ by a vector of norm less than $nl(Q)$, so Lemma
$\ref{lem3}$ implies that $|D(\widetilde{Q}) - D(\widetilde{Q'})|
\leq 3n^{3}\omega(l(Q))$. Hence, to show $\eqref{eq41}$ it is enough
to prove that
\begin{equation}\label{eq42} \lim_{|Q| \rightarrow 0} \frac{|A \cap
\phi^{-1}(Q)| - |A \cap \widetilde{Q}|}{|Q|} = 0
\end{equation} We study first the case $n=2$. We are going to decompose
$\phi^{-1}(Q)$ into squares as follows. Let $Q_{0}$ be the maximal
square with sides parallel to the axis contained in $\phi^{-1}(Q)$
and write $\mathfrak{F}_{0}=\left\{Q_{0}\right\}$. Observe that the
ratio of the area of $Q_{0}$ to that of $\phi^{-1}(Q)$ is
$C=1/(1+\sin(2\alpha))$ and that $0.5 \leq C \leq 4- 2 \sqrt{3}$.
Then $\phi^{-1}(Q) \backslash Q_{0}$ is the union of eight
right-angled triangles whose hypothenuse is contained in $\partial
\phi^{-1}(Q)$. See Figure 1. \begin{figure}[ht]
\begin{center}
\epsfig{figure=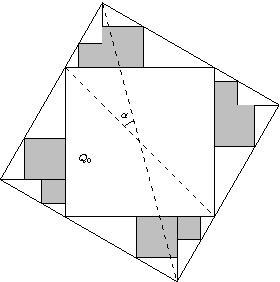} \caption{The shaded square are the eight
elements of $\mathfrak{F}_{1}$ for $\alpha=\pi/6$.}
\end{center}
\end{figure}
Take again the maximal square with sides parallel to the axis
contained in each triangle, obtaining a family $\mathfrak{F}_{1}$ of
eight squares of total area $(1-C)^{2}|Q|$. Thus $\phi^{-1}(Q)
\backslash (\mathfrak{F}_{0} \cup \mathfrak{F}_{1})$ is the union of
16 right-angled triangles and we continue inductively, constructing,
for $k=1,2...$, a family $\mathfrak{F}_{k}$ of $2^{k+2}$ squares of
total area $C^{k-1}(1-C)^{2}|Q|$. Observe that
\[D(\phi^{-1}(Q)) - D(\widetilde{Q}) = \sum_{k} \sum_{R \in
\mathfrak{F}_{k}} \frac{|R|}{|Q|} (D(R)- D(\widetilde{Q}))\] Let $R$
be a square in $\mathfrak{F}_{k}$. Since $A$ is smooth, there exists
a constant $C_{1}>0$ such that $|D(R)-D(\widetilde{Q})| \leq C_{1} k
\omega(l(Q))$. We deduce that
\begin{align}|D(\phi^{-1}(Q)) - D(\widetilde{Q})| &\leq C_{1}
\omega(l(Q)) \sum_{k} k \frac{\left|\bigcup_{\mathfrak{F}_{k}} R
\right|}{\left|\phi^{-1}(Q)\right|} \leq \nonumber \\ & \leq
C_{1}(1-C)^{2}C^{-1}\omega(l(Q)) \sum_{k} k C^{k} = C_{1}
\omega(l(Q)) \nonumber
\end{align} This
implies $\eqref{eq42}$ and finishes the proof in dimension $2$ when
$\phi$ is a rotation.

We now study the higher dimensional case $n>2$. Recall that
$\phi^{-1}$ is a rotation on a two dimensional subspace $E$ and
$\phi^{-1}$ is the identity on its orthogonal complement. Without
loss of generality, we may assume that $E$ is generated by the two
first vectors of the canonical basis of $\R^{n}$. Consider the
orthogonal projection $\Pi$ of $\R^{n}$ onto $E$ and decompose
$\Pi(\phi^{-1}(Q))$ as in the two dimensional case, that is
$\Pi(\phi^{-1}(Q)) =
\displaystyle{\bigcup_{k=0}^{\infty}\mathfrak{F}_{k}}$, where
$\mathfrak{F}_{k}$ is, as before, the union of $2^{k+2}$ (2
dimensional) squares with sides parallel to the axis of total area
$C^{k-1}(1-C)^{2}l(Q)^{2}$. Since $\phi^{-1}(Q)= \Pi(\phi^{-1}(Q))
\times B$ where $B$ is a cube in $\R^{n-2}$ with sides parallel to
the axis, we have $\phi^{-1}(Q)=
\displaystyle{\bigcup_{k=0}^{\infty} G_{k}}$ where
$G_{k}=\displaystyle{\bigcup_{R \in \mathfrak{F}_{k}} R \times B}$.
Using the smoothness condition one can show that there exists a
constant $C_{2}>0$ such that $|D(R \times B) - D(\widetilde{Q})|
\leq C_{2} k \omega (l(Q))$ for any square $R \in \mathfrak{F}_{k}$.
Then \[|D(\phi^{-1}(Q)) - D(\widetilde{Q})| \leq C_{2}\omega(l(Q))
\sum_{k=0}^{\infty} k \frac{|G_{k}|}{|Q|}\] Since $|G_{k}|=
(1-C)^{2}C^{k-1} |Q|$ with $0.5 < C < 4 - 2 \sqrt{3}$ we deduce that
$|D(\phi^{-1}(Q))-D(\widetilde{Q})| \leq C_{2} \omega (l(Q))$ and
thus $\eqref{eq41}$ is satisfied, completing the proof in the case
that $\phi$ is a rotation.

Let us now assume that $\phi$ is a dilation, that is, $\phi$ has a
diagonal matrix. Without loss of generality we can assume
$\phi^{-1}(x_{1}, ... , x_{n}) = (\lambda x_{1}, x_{2},...,x_{n})$
for some $\lambda \in \R$. Assume $\lambda > 1$. To prove
$\eqref{eq41}$ it is sufficient to find a constant $C(\lambda)>0$
such that for any cube $Q \subset \R^{n}$ and any cube
$\widetilde{Q} \subset \phi^{-1}(Q)$ with $l(\widetilde{Q})=l(Q)$ we
have
\begin{equation}\label{eq43} |D(\phi^{-1}(Q)) - D(\widetilde{Q})|
\leq C(\lambda) \omega(l(Q))
\end{equation} The proof of $\eqref{eq43}$ resembles that of part (a) of Lemma
$\ref{lem3}$. We can assume that $\widetilde{Q}$ is the unit cube.
Let $[\lambda]$ be the integer part of $\lambda$ and write the
interval $[[\lambda],\lambda)$ as a union of maximal dyadic
intervals $\{I_{k}\}$ with $|I_{k}|=2^{-k}$, that is,
$[[\lambda],\lambda) = \bigcup I_{k}$. Consider $R_{k}= I_{k} \times
[0,1]^{n-1}, \quad k=1,2,...$. Observe that \[|D(\phi^{-1}(Q)) -
D(\widetilde{Q})| = \sum_{j=0}^{[\lambda]-1} \frac{1}{\lambda}
\left(D([j,j+1) \times [0,1]^{n-1}) - D(\widetilde{Q})\right) +
\sum_{k=1}^{\infty} \frac{2^{-k}}{\lambda} (D(R_{k}) -
D(\widetilde{Q}))\] For $j=0,1,..., [\lambda]-1$, we have
$|D([j,j+1) \times [0,1]^{n-1}) - D(\widetilde{Q})| \leq \lambda
\omega(1)$. Since $R_{k}$ can be split into a family of dyadic cubes
of generation $k$, we deduce that $|D(R_{k}) - D(\widetilde{Q})|
\leq (\lambda + 1 +k)\omega(1)$. Therefore
\[|D(\phi^{-1}(Q))-D(\widetilde{Q})| \leq (\lambda + 1 + 3/\lambda)
\omega(1)\] which proves $\eqref{eq43}$. An analogous argument can
be used in the case $\lambda<1$.
\end{proof}

\begin{remark}\label{rem1}
The first part of the proof shows that there exists a constant
$C=C(n) >0$ such that for any rotation $\phi$ in $\R^{n}$ and any
$\omega$-smooth set $A \subset \R^{n}$, its image $\phi(A)$ is
$C\omega$-smooth. When $\phi$ is a dilation in a single direction
with parameter $\lambda \in \R$ and $A \subset \R^{n}$ is an
$\omega$-smooth set, the proof shows that $\phi(A)$ is
$C(\lambda)\omega$-smooth, with $C(\lambda) \leq 4 (\lambda +
1/\lambda)$.
\end{remark}
\begin{remark}\label{rem2}
Let $\left\{T_{i}\right\}$ be a countable family of linear
isomorphisms in $\R^{n}$ for which there exists a constant $M
> 0$ such that $M^{-1} ||x|| \leq ||T_{i}(x)|| \leq M ||x||$ for any
$x \in \R^{n}$ and $i=1,2...$ Then there exists a constant
$C=C(M,n)>0$ such that for any $\omega$-smooth set $A$ and any $i$
one has \[\frac{||T_{i}(Q) \cap A| - |T_{i}(Q') \cap A||}{|Q|} \leq
C\omega(l(Q))\]
\end{remark}
\begin{remark}\label{rem3}
Proposition $\ref{prop41}$ and part (a) of Lemma $\ref{lem3}$ give
that affine mappings preserve smooth sets.
\end{remark}

We could have defined smooth sets using the grid of dyadic cubes or,
in the opposite direction, using the grid of all cubes, even without
taking them parallel to the axis. The previous results imply that
both grids would lead to equivalent definitions.

\begin{corollary}
Let $A$ be a measurable set in $\R^{n}$. The following are
equivalent:
\begin{description}
\item [(a)] For any $\varepsilon > 0$, there exists $\delta >0$ such
that $|D(Q) - D(Q')| \leq \varepsilon$ for any pair of consecutive
dyadic cubes $Q, Q'$, of the same side length $l(Q)=l(Q') < \delta$.
\item [(b)] For any $\varepsilon > 0$, there exists $\delta >0$ such
that $|D(Q) - D(Q')| \leq \varepsilon$ for any pair of consecutive
cubes $Q, Q'$ with sides non necessarily parallel to the axis, of
the same side length $l(Q)=l(Q') < \delta$.
\item [(c)] $A$ is a smooth set.
\end{description}
\end{corollary}

\begin{proof}
Since any cube in $\R^{n}$ is the affine image of a dyadic cube,
Lemma $\ref{lem3}$ shows that (a) implies (b). The other
implications are obvious.
\end{proof}

Observe that bilipschitz mappings do not preserve smoothness in
general. Applying locally Proposition $\ref{prop41}$, we can extend
it to certain diffeomorphisms, but we need extra assumptions to
guarantee that the local bounds that we obtain are satisfied
uniformly. We are now ready to proceed with the proof of Theorem
$\ref{th3}$.

\begin{proof}[of Theorem $\ref{th3}$]
Since $\phi$ is bilipschitz, $|\phi(Q)|$ is comparable to $|Q|$.
Also, $|J\phi|$ is uniformly bounded from above and below.
Therefore, $J\phi^{-1}$ is uniformly continuous as well. We will
also need that
\begin{equation}\label{eq44}
\lim_{|Q| \rightarrow 0} \frac{|\phi(Q)| - |\phi(Q')|}{|Q|} = 0
\end{equation} To show this, observe that this quantity is \[\frac{1}{|Q|}
\left(\int_{Q} J\phi - \int_{Q'} J\phi \right)\] which tends to $0$
uniformly when $l(Q) \rightarrow 0$ because of the uniform
continuity of $J\phi$.

We first show that (b) is equivalent to (c). A change of variables
gives that \[|\phi^{-1}(A) \cap Q| - |\phi^{-1}(A) \cap Q'| = \int
J\phi^{-1}(x) \left(\1_{A \cap \phi(Q)}(x) - \1_{A \cap \phi(Q')}
(x) \right) dx\] Let $p(Q)$ be a point in $\overline{\phi(Q)} \cap
\overline{\phi(Q')}$. Given $\varepsilon
> 0$, if $l(Q)$ is sufficiently small one has $||J\phi^{-1}(x) -
J\phi^{-1}(p(Q))|| < \varepsilon$ for any $x \in \phi(Q)$. Hence,
the uniform continuity of $J\phi^{-1}$ gives us that \[\lim_{|Q|
\rightarrow 0} \frac{|\phi^{-1}(A) \cap Q| - |\phi^{-1}(A) \cap
Q'|}{|Q|} = \lim_{|Q| \rightarrow 0} \frac{\left(|A \cap \phi(Q)| -
|A \cap \phi(Q')|\right) J\phi^{-1}(p(Q))}{|Q|}\] Let $D(\phi(Q))$
be the density of $A$ in $\phi(Q)$, that is, $D(\phi(Q))= |A \cap
\phi(Q)|/ |\phi(Q)|$. Applying $\eqref{eq44}$ we have \[\lim_{|Q|
\rightarrow 0} \frac{|\phi^{-1}(A) \cap Q| - |\phi^{-1}(A) \cap
Q'|}{|Q|} = \lim_{|Q| \rightarrow 0} \left(D(\phi(Q)) - D(\phi(Q'))
\right) J\phi^{-1}(p(Q))\] Since $J\phi^{-1}$ is uniformly bounded
both from above and below, we deduce that (b) and (c) are
equivalent.

We now show that (a) implies (c). Observe that, applying
$\eqref{eq44}$, it is sufficient to show
\begin{equation}\label{eq45} \lim_{|Q| \rightarrow 0} \frac{|A \cap
\phi(Q)| - |A \cap \phi(Q')|}{|Q|} = 0
\end{equation} Let $z(Q)$ be a point in $Q$ with dyadic coordinates. Let $T=T(Q)$ be
the affine mapping defined by $T(x) = \phi(z(Q)) + D\phi (z(Q))
(x-z(Q))$, for any $x \in \R^{n}$, where $D\phi$ denotes the
differential of $\phi$. Given $\varepsilon > 0$, the uniform
continuity of $J\phi$ tells that $|\phi(x) - T(x)| \leq \varepsilon
l(Q)$, for any $x \in Q \cup Q'$ if $l(Q)$ is sufficiently small.
There thus exists a constant $C_{1}(n)>0$ such that if $l(Q)$ is
sufficiently small then \[|\left(\phi(Q) \backslash T(Q)\right) \cup
\left(T(Q) \backslash \phi(Q)\right)| \leq C_{1}(n) \varepsilon
|Q|\] and similarly for $Q'$. So we deduce that $\eqref{eq45}$ is
equivalent to
\begin{equation}\label{eqlast} \lim_{|Q| \rightarrow 0} \frac{|A \cap T(Q)| - |A \cap
T(Q')|}{|Q|} = 0\end{equation} Now $\eqref{eqlast}$ follows from
Remark $\ref{rem2}$ because, since $\phi$ is bilipschitz, there
exists a constant $M>0$ such that $M^{-1}||x|| \leq
||D\phi(z(Q))(x))|| \leq M ||x||$ for any $x \in \R^{n}$ and any
cube $Q$ in $\R^{n}$. This finishes the proof that (a) implies (c).
The proof that (b) implies (a) follows applying the previous part to
$\phi^{-1}$.
\end{proof}

$\quad$

$\quad$

Artur Nicolau and Daniel Seco

Departament de Matem{\`a}tiques

Universitat Aut{\`o}noma de Barcelona

08193 Bellaterra

Spain

artur@mat.uab.cat

dseco@mat.uab.cat
\end{document}